\documentclass{article}
\usepackage{amssymb}
\usepackage{amsmath}
\usepackage{enumitem}
\usepackage{tikz-cd}
\usepackage{hyperref}

\usepackage[titletoc]{appendix}
\usepackage{amsthm}
\usepackage{parskip}
\usepackage{amsfonts}
\usepackage{caption}
\usepackage{mathrsfs}
\usepackage{afterpage}
\usepackage{float}
\usepackage{multirow}
\usepackage{graphicx}
\usepackage{color}
\usepackage[margin=1.0in]{geometry}
\newcommand*{\SHom}{\mathscr{H}\kern -.5pt om}
\newcommand*{\SExt}{\mathscr{E}\kern -.5pt xt}

\begin{document}
\title{On the exponential type conjecture}
\date{}
\author{Zihong Chen}
\maketitle
\theoremstyle{definition}
\newtheorem{mydef}{Definition}[section]
\numberwithin{mydef}{section}
\newtheorem{rmk}[mydef]{Remark}
\newtheorem{conj}[mydef]{Conjecture}
%\numberwithin{rmk}{subsection}
\theoremstyle{plain}
\newtheorem{cor}[mydef]{Corollary}
%\numberwithin{cor}{subsection}
\newtheorem{lemma}[mydef]{Lemma}
%\numberwithin{lemma}{subsection}
\newtheorem{thm}[mydef]{Theorem}
%\numberwithin{thm}{subsection}
\newtheorem{prop}[mydef]{Proposition}
%\numberwithin{prop}{subsection}

\makeatletter
\def\thm@space@setup{%
  \thm@preskip=\parskip \thm@postskip=0pt
}
\makeatother

\begin{abstract}
We prove that the small quantum $t$-connection on a closed monotone symplectic manifold is of exponential type and has quasi-unipotent regularized monodromies at $t=0$. This answers a conjecture of Katzarkov-Kontsevich-Pantev and Galkin-Golyshev-Iritani for those classes of symplectic manifolds. The proof follows a reduction to positive characteristics argument, and the main tools of the proof are Katz's local monodromy theorem in differential equations and quantum Steenrod operations in equivariant Gromov-Witten theory with mod $p$ coefficients.
\end{abstract}

\renewcommand{\theequation}{1.\arabic{equation}}
\setcounter{equation}{0}
\section{Introduction}
Fix a closed (i.e. compact without boundary) symplectic manifold $(X,\omega)$ of dimension $2n$ that is \emph{monotone}, i.e.
\begin{equation}
[\omega]=c_1(X)\in H^2(X;\mathbb{R}).   
\end{equation}
Unless otherwise specified, by a graded module/vector space/algebra we will mean a $\mathbb{Z}/2$-graded module/vector space/algebra; accordingly, when we say an element has degree $d$, it should be understood as having degree $d$ mod $2$.\par\indent
Let $t$ be a formal variable of degree $2$. This paper is concerned with the following version of a single variable quantum connection with quadratic pole introduced by Dubrovin, cf.  \cite{Dub1} and \cite{Dub2}.
\begin{mydef}
The \emph{quantum $t$-connection} is the formal connection on $H^*(X;\mathbb{C})((t))$ defined by
\begin{equation}
\nabla^{QH}_{\frac{\partial}{\partial t}}:=\frac{\partial}{\partial t}+\frac{\mu}{t}-\frac{c_1\star}{t^2},    
\end{equation}
where $\mu$ is the grading operator defined by $\mu|_{H^k(X)}:=\frac{k-n}{2}$ and $\star: H^*(M,\mathbb{C})\otimes H^*(M,\mathbb{C})\rightarrow H^*(M,\mathbb{C})$ denotes the small quantum cup product. 
\end{mydef}
At $t=\infty$, (1.2) has a regular singularity, making the local structure relatively simple; for instance, there exists a basis of local flat sections extending any basis of $H^*(X;\mathbb{C})$ at $t=\infty$. On the other hand, the local structure of (1.2) at $t=0$ is much more complicated due to its apparent second order pole, which in general cannot be reduced to a simple pole by gauge transformations. \par\indent
There is a classical and rich theory on the formal classification of connection with (irregular) singularity, cf. \cite{BV}, \cite{Lev}, \cite{Mal}. Restricting to the case of a formal connection with quadratic pole, the \emph{Hukuhara-Levelt-Turittin theorem} (cf. \cite[Section 3]{Lev},\cite[Theorem 2.1]{Mal}) implies that any such connection $(M,\nabla)$ over $\mathbb{C}((t))$ admits a finite decomposition, possibly after a finite base change $p_N(t)=t^N: \mathbb{C}((t))\rightarrow \mathbb{C}((t))$, of the form
\begin{equation}
    p_N^*(M,\nabla)\cong\bigoplus_{\lambda}(\mathbb{C}((t)),\frac{\partial}{\partial t}-\frac{\lambda}{t^2})\otimes \nabla^{reg,\lambda}_{\frac{\partial}{\partial t}},
    \end{equation}
where $\nabla^{reg,\lambda}_{\frac{\partial}{\partial t}}$ is a connection with regular singularity. Moreover, the integer $N$ (called the \emph{ramification index}) is an invariant of the connection and such a decomposition is unique up to isomorphism.\par\indent
The following conjecture, which distinguishes the formal structure of the quantum connection among general connections with quadratic singularities, first appeared in \cite[Conjecture 3.4]{KKP} and later in a different context in \cite[Section 2.5]{GGI}.
 \begin{conj}[The exponential type conjecture]
    The quantum $t$-connection $(H^*(X;\mathbb{C})((t)),\nabla^{QH}_{\frac{\partial}{\partial t}})$ of (1.2) admits a finite direct sum decomposition (1.3) with ramification index $1$. Moreover, the monodromies of the regular parts $\nabla^{reg,\lambda}$ are given by roots of unity. In fancier language, this says that the quantum connection has \emph{unramified exponential type} and \emph{quasi-unipotent regularized monodromies}.
\end{conj}
 Conjecture 1.2 is a key first step in setting up a (so far largely conjectural) nice theory of Hodge structures for quantum cohomology on a general closed monotone symplectic manifold. One main consequence of being of exponential type is that locally near $t=0$ analytic solutions to the quantum connection exist in a large range of directions (known as a \emph{Stokes sector}). Studying the interactions of (the analytic continuation) of these solutions with the solutions near $t=\infty$, as well as comparisons of solutions across different Stokes sectors, lie at the heart of deep questions in Gromov-Witten theory such as the the Dubrovin conjecture and the Gamma conjectures. These topics are outside the scope of this paper, and we refer the readers to \cite{CDG}, \cite{GGI}, \cite{SS} for a detailed discussion. \par\indent
The main result of this paper is the following.
\begin{thm}
Conjecture 1.2 holds for all closed monotone symplectic manifolds $(X,\omega)$.    
\end{thm}
\begin{rmk}
\begin{enumerate}[label=\roman*)]
 \item The first substantial progress towards this conjecture was made in a series of papers of Pomerleano-Seidel \cite{PS1},\cite{PS2}, where they proved Conjecture 1.2 assuming the existence of a smooth anticanonical divisor whose complement is Weinstein. Their approach involved passing to the Fukaya category of the divisor complement combined with prior results on the non-commutative Gauss-Manin connection \cite{PVV}.
    \item Independently, Conjecture 1.2 was later proved by the author, cf. \cite[Theorem 1.2]{Che}, conditional on the assumption that $X$ satisfies Abouzaid's generation criterion \cite{Abo}. This assumption was used in loc.cit. to study quantum Steenrod operations via the monotone Fukaya category and to prove a key compatibility relation with the elementary splitting of the quantum connection, cf. Corollary 1.10 of loc.cit.. The current paper is built on the approach of \cite{Che}, but completely bypasses Fukaya category and vastly strengthens and simplifies the argument.

    \item The approach of this paper also provides an algorithm for finding a basis that puts the quantum $t$-connection into canonical form, i.e. for which the associated connection matrix of $\nabla^{QH}_{\frac{\partial}{\partial t}}$ is block diagonal where each block is of the form
    \begin{equation}
     \frac{-\lambda I_n}{t^2}+\frac{A_{-1}}{t}+A_0+A_1t+\cdots,\;\lambda\in\mathbb{C},A_i\in \mathbf{M}_n(\mathbb{C}).
    \end{equation}
    Namely, one first runs the splitting argument of Lemma 2.1, which involves explicitly solving for a power series by induction. Then, for each connection arising as a summand of the splitting, the proof of Katz's local monodromy theorem (Theorem A.1) guarantees that any basis generated by a cyclic vector (which can be find explicitly by an algorithm of Katz \cite{Ka2}) puts the corresponding summand into canonical form. 
    \item Finally, we remark that while the result of this paper subsumes various cases of the exponential type conjecture proved in \cite{Che} or \cite{PS1},\cite{PS2}, the approach of Pomerleano-Seidel leads to certain quantitative conclusions (such as bounds on the sizes of Jordan blocks of the regularized monodromies) that are unique to their approach. To the author's awareness, a further study of these invariants and a comparison of the aforementioned approaches will appear in an announced work of Pomerleano-Seidel. 
\end{enumerate}
   
\end{rmk}  
The proof of Theorem 1.3 follows a reduction mod $p$ argument, by bringing together Katz's local monodromy theorem in the theory of differential equations and quantum Steenrod operations in equivariant Gromov-Witten theory with mod $p$ coefficients. Below we give a very brief introduction to these subjects.\par\indent
\textbf{$p$-curvature and Katz's local monodromy theorem}. To introduce the notion of $p$-curvature, we first make the following simple observation. Let $X$ be a (formal) scheme over a field $\mathbf{k}$ of characteristic $p>0$, and $D: \mathcal{O}_X\rightarrow \mathcal{O}_X$ be a (continuous) $\mathbf{k}$-linear derivation on $X$. Then its $p$-th iteration 
\begin{equation}
D^p:=\overbrace{D\circ D\circ\cdots\circ D}^{p\;times}: \mathcal{O}_X\rightarrow \mathcal{O}_X  
\end{equation}
remains a derivation (an easy application of the Leibniz rule), sometimes called the \emph{$p$-th power of $D$}. \par\indent
In the above setting, further fix $(M,\nabla)$ an $\mathcal{O}_X$-module with connection. Just as curvature measures the incompatibility of a connection with Lie brackets, $p$-curvature measures the incompatibility of a connection (in characteristic $p$) with $p$-th powers. 
\begin{mydef}
The \emph{$p$-curvature} of $\nabla$ along $D$ is defined to be
\begin{equation}
F^{\nabla}_D:=\nabla^p_{D}-\nabla_{D^p}: M\rightarrow M.
\end{equation}    
\end{mydef}
Just like curvature, it turns out that $F^{\nabla}_D$ is $\mathcal{O}_X$-linear, and thus (1.6) can be packaged into a map 
\begin{equation}
F^{\nabla}: \mathrm{Der}_{\mathbf{k}}(\mathcal{O}_X)\rightarrow \mathrm{End}_{\mathcal{O}_X}(M). 
\end{equation}
Moreover, (1.7) is additive and Frobenius $p$-linear (i.e. $F^{\nabla}_{g D}=g^pF^{\nabla}_D$). We refer the readers to \cite[Section 5]{Ka1} for proofs of these properties and a more detailed introduction.  \par\indent
One of the most important features of $p$-curvature is that it gives a criterion (known as \emph{Cartier descent}, cf. \cite[Theorem 5.1]{Ka1}) for when a differential equation in characteristic $p$ has a full set of algebraic solutions. Namely, this is the case if and only if its $p$-curvature vanishes. This starkly contrasts the case of characteristic $0$, where no such simple criterion for the existence of algebraic solutions to a general differential equation (with singularity) is known. A well-known (and still wildly open) conjecture of Grothendieck and Katz (cf. \cite{Ka3}, \cite{Ka4}) predicts that the existence of algebraic solutions to a differential equation in characteristic $0$ is encoded in the $p$-curvature of its mod $p$ reduction across all primes $p$. \par\indent
The relevant ingredient we need from the arithmetic theory of differential equations is  a result of Katz \cite[Theorem 13.0]{Ka1} known as the \emph{local monodromy theorem}, which can be viewed as (a positive answer to) the local version of the Grothendieck-Katz conjecture, cf. Appendix A for more details. More specifically, it gives a criterion for when a connection is regular singular in terms of the $p$-curvatures of its reduction mod $p$'s. We will see that a suitable application of Katz's theorem reduces Conjecture 1.2 to studying the $p$-curvature of the quantum connection. \par\indent
\textbf{Quantum Steenrod operations}. The second main ingredient of the proof is the collection of \emph{quantum Steenrod operations}. First introduced by Fukaya \cite{Fuk}, and systematically developed by Wilkins \cite{Wil}, these operations quantize classical Steenrod operations on the mod $p$ cohomology of a symplectic manifold via equivariant Gromov-Witten theory with mod $p$ coefficients. They have seen various links and applications to Hamiltonian dynamics \cite{Shel}\cite{CGG}\cite{Rez}, arithmetic mirror symmetry \cite{Sei3} and representation theory \cite{Lee} \cite{BL}.\par\indent
An important feature of quantum Steenrod operations is their relation with the quantum connection. For instance, Seidel-Wilkins \cite{SW} showed that they are covariantly constant with respect to the quantum connection. Work of Lee \cite{Lee} further showed that quantum Steenrod operations of divisor classes agree with the $p$-curvature of the quantum connection in the case of symplectic resolutions, and envisioned this to be a more general phenomenon. These properties give a geometric perspective on the $p$-curvature of the quantum connection, and play a crucial role in the proof of this paper. We give a more detailed review of quantum Steenrod operations in Section 3.\par\indent
We now sketch the main steps of the argument.
\begin{enumerate}[label=\arabic*)]
    \item By the elementary splitting lemma, cf. Lemma 2.1, there exists an $R\subset \overline{\mathbb{Q}}$ finitely generated over $\mathbb{Z}$ together with a decomposition 
   \begin{equation}
H^*(X;R)[[t]]=\bigoplus_{\lambda\in\mathrm{spec}(c_1\star)} H^*(X;R)[[t]]_{-\lambda},
\end{equation}
of the quantum $t$-connection, where the quadratic pole term $-c_1\star$ of $\nabla^{QH}_{\frac{\partial}{\partial t}}$ restricted to the component $H^*(X;R)[[t]]_{-\lambda}$ has a single eigenvalue $-\lambda$. We remark that the relevant properties of $R$ for applying the `spreading out' argument in Katz's theorem outlined in 2) are that $R$ has infinitely many maximal ideals the intersection of which is zero, and that each residue field is finite (such an $R$ is non-unique, and any choice works).  
     \item Applying (a formal version of) Katz's local monodromy theorem, cf. Theorem A.1, one can reduce the proof of Theorem 1.3 to showing that 
     \begin{equation}
      F_{t^2\frac{\partial}{\partial t}}^{QH,-\lambda}+\lambda^p   
     \end{equation}
     is a nilpotent endomorphism of $H^*(X,R/\mathfrak{m})[[t]]_{-\lambda}$  for all maximal ideals $\mathfrak{m}\subset R$. In (1.9), $F^{QH,-\lambda}_{t^2\frac{\partial}{\partial t}}$ denotes the \emph{$p$-curvature} (cf. Definition 1.5) of the connection $\nabla^{QH}|_{H^*(X;R)[[t]]_{-\lambda}}$ base changed to $R/\mathfrak{m}$, where $p$ denotes the characteristic of the (finite) field $R/\mathfrak{m}$.
     \item A key insight, first observed by Jae Hee Lee when studying quantum connections of symplectic resolutions \cite{Lee} and later adapted to the setting of closed monotone symplectic manifolds by Paul Seidel (the author learned of the latter result via private communication), is that the operator 
     \begin{equation}
       Q\Sigma_{c_1}+F^{QH}_{t^2\frac{\partial}{\partial t}}\in \mathrm{End}(H^*(X,R/\mathfrak{m})[[t,\theta]])
       \end{equation}
     is nilpotent, where $\theta$ is an odd square zero formal variable commuting with $t$ (more precisely, it should be interpreted as the generator of $H^1(B\mathbb{Z}/p;\mathbb{F}_p)$). In (1.10), $Q\Sigma$ denotes the \emph{quantum Steenrod operations}, cf. \cite{Fuk}, \cite{Wil} and \cite{SW}, which is a Frobenius $p$-linear algebra action
     \begin{equation}
      Q\Sigma: QH^*(X,R/\mathfrak{m})\otimes H^*(X,R/\mathfrak{m})[[t,\theta]]\rightarrow H^*(X,R/\mathfrak{m})[[t,\theta]]
     \end{equation}
     constructed out of $\mathbb{Z}/p$-Borel-equivariant Gromov-Witten invariants with mod $p$ coefficients, cf. Section 3 for more details. Since $Q\Sigma$ commutes with the quantum connection (cf. Theorem 3.2, Lemma B.2) and hence its $p$-curvature, we are further reduced to showing that $Q\Sigma_{c_1}-\lambda^p$ is a nilpotent operation when restricted to $H^*(X,R/\mathfrak{m})[[t,\theta]]_{-\lambda}$. 
     \item In the final step, we conclude the proof by establishing a compatibility relation between $Q\Sigma$ and the decomposition (1.8), after base changed to $R/\mathfrak{m}$, cf. Corollary 4.4. We remark that this statement is essentially \cite[Corollary 1.10]{Che}, which was proved in loc.cit. using Fukaya categories and cyclic open-closed maps under the assumption that $X$ satisfies Abouzaid's generation criterion. In this paper, we prove it for all closed monotone symplectic manifolds using purely closed-string methods. 
     
\end{enumerate}
In Appendix C, we demonstrate that a similar strategy can be applied to give a short proof that the canonical connection associated with the derived category of an isolated singularity (which often arises as the mirror B-side of a closed monotone symplectic manifold) has exponential type, recovering a simple case of a classical result of Sabbah \cite{Sab}. 

\subsection*{Acknowledgements}
First and foremost, I would like to thank my PhD advisor Paul Seidel for his patient and helpful guidance throughout my graduate study, as well as numerous enlightening conversations regarding this paper. I would like to thank Kai Hugtenburg for answering my various questions regarding formal decomposition of connections; and Jae Hee Lee for his insightful work that brought together $p$-curvature of quantum connections and quantum Steenrod operations, an idea that enabled the approach in this paper. I would also like to thank the anonymous referees for their invaluable comments on the manuscript, which has greatly improved its exposition and readability. This research was partially supported by the Simons Foundation, through a Simons Investigator grant (256290).

\renewcommand{\theequation}{2.\arabic{equation}}
\setcounter{equation}{0}
\section{Formal connections with quadratic singularities}
This section collects some basic facts about general connections with quadratic singularities on the formal affine line. The main output is Proposition 2.5 (which to the author's awareness is new), which gives a characterization of the summands in Levelt's splitting (Lemma 2.1) in terms of covariantly constant endomorphisms. \par\indent  
Let $R$ be an integral domain and $V$ be a finite rank free $R$-module. Let $t$ be a formal variable. A central object of study in this paper is a formal connection on $V((t))$ with a quadratic singularity at $t=0$
\begin{equation}
\nabla_\frac{\partial}{\partial t}=\frac{\partial}{\partial t}+\frac{A_0}{t^2}+\frac{A_1}{t}+\cdots: V((t))\rightarrow V((t)),
\end{equation}
where $A_i\in \mathrm{End}(V)$. An example is the quantum $t$-connection of (1.2). \emph{Throughout this section, we assume that $A_0\in \mathrm{End}(V)$ has a Jordan decomposition over $R$}, i.e. 
\begin{itemize}
    \item The characteristic polynomial $p_{A_0}(t)$ of $A_0$ splits completely over $R$, i.e. $p_{A_0}(t)=\prod_{i=1}^m(t-\lambda_i)^{k_i}$, where $\lambda_i\in R$. 
    \item Moreover, $V$ decomposes into a direct sum   
    \begin{equation}
     V=\bigoplus_i V_{\lambda_i},   
    \end{equation}
where each $\lambda_i$-generalized eigenspace $V_{\lambda_i}:=\mathrm{ker}(A_0-\lambda_iI)^{k_i}$ is a finite rank free $R$-submodule of $V$. 
\end{itemize}
Lemma 2.1-Corollary 2.4 below are well-known, cf. for instance \cite[Section 2]{Lev}, \cite[Lemma 2.13, Lemma 2.14, Lemma 2.16]{Hug}. However, they are usually stated and proved over $\mathbb{C}$, and we include slight modifications of the proofs to adapt to the more general setting suitable for our needs. 
\begin{lemma} (Elementary splitting lemma)
Assume that the difference between two distinct eigenvalues of $A_0$ is invertible in $R$. Write $\{e_j\}$ for a generalized eigenbasis of $A_0$, and write $\nabla=\frac{\partial}{\partial t}+\frac{A_0}{t^2}+\frac{A_1}{t}+\cdots$ in this basis. In particular, $A_0$ is in Jordan normal form with diagonal entries the eigenvalues of $A_0$. Then, there exists a basis $\{v_j\}$ of $V[[t]]$ such that 
\begin{enumerate}[label=\arabic*)]
    \item  $v_j|_{t=0}=e_j$.
    \item Write $\nabla=\frac{\partial}{\partial t}+\frac{1}{t^2}\sum_{i=0}^{\infty} \tilde{A}_i t^i$ in the basis $\{v_j\}$, then $\tilde{A}_0=A_0$ and all $\tilde{A}_i$'s respect the generalized eigen-decomposition of $A_0$.
\end{enumerate}
\end{lemma}
\emph{Proof}. The proof is by induction. Suppose that $A_1,\cdots,A_{r-1}$ all respect the generalized eigen-decomposition of $A_0$. We make a formal gauge transformation by
\begin{equation}
P=\mathrm{id}+t^rT_r,    
\end{equation}
for some undetermined $T_r\in \mathrm{End}(V)$. Then, the gauged transformed connection matrix is
\begin{equation}
P^{-1}(\sum_{i\geq 0}A_it^{i-2})P-P^{-1}\frac{dP}{dt}=\sum_{0\leq i\leq r-1}A_it^{i-2}+(A_r+[A_0,T_r])t^{r-2}+O(t^{r-1}). 
\end{equation}
Consider the endomorphism of $\mathrm{End}(V)$ given by $X\mapsto [X,A_0]$, which respects the decomposition 
\begin{equation}
\mathrm{End}(V)=\bigoplus_{\lambda,\lambda'}\mathrm{Hom}(V_{\lambda},V_{\lambda'}).     
\end{equation}
The $(\lambda,\lambda')$-component of this endomorphism is given by
\begin{equation}
X^{\lambda,\lambda'}\mapsto X^{\lambda,\lambda'}A_0^{\lambda'}-A_0^{\lambda}X^{\lambda,\lambda'}, 
\end{equation}
where $A_0^{\lambda}$ denotes the restriction of $A_0$ to $V_{\lambda}$. Rewrite (2.6) as
\begin{equation}
X^{\lambda,\lambda'}\mapsto (\lambda-\lambda')X^{\lambda,\lambda'}+\big(X^{\lambda,\lambda'}(A^{\lambda'}_0-\lambda')-(A_0^{\lambda}-\lambda)X^{\lambda,\lambda'}\big).
\end{equation}
If $\lambda\neq \lambda'$, then by assumption $\lambda-\lambda'$ is invertible and moreover $A_0^{\lambda}-\lambda, A_0^{\lambda'}-\lambda'$ are nilpotent by definition. This implies that (2.6) is an isomorphism. \par\indent
In particular, there is a unique $T_r$ such that
\begin{enumerate}[label=\arabic*)]
    \item $T_r$ is off-block-diagonal with respect to the generalized eigendecomposition of $A_0$ and
    \item the off-block-diagonal terms of $-[A_0,T_r]$ agrees with those of $A_r$. 
\end{enumerate}
Thus, applying the gauge transformation (2.3) puts the new $\tilde{A}_r=A_r+[A_0,T_r]$ into the desired form. Now, we can proceed with the induction and note that the infinite product of all the gauge transformations involved
\begin{equation}
 (\mathrm{id}+T_1t)(\mathrm{id}+T_2t^2)(\mathrm{id}+T_3t^3)\cdots  
\end{equation}
is a well-defined formal gauge transformation. Applying (2.8) to the basis $\{e_j\}$ gives the desired basis $\{v_j\}$. \qed\par\indent
In particular, Lemma 2.1 implies that for any formal connection $\nabla_{\frac{\partial}{\partial t}}$ with a quadratic pole, there exists a decomposition
\begin{equation}
V[[t]]=\bigoplus_{\lambda} V[[t]]_{\lambda}   
\end{equation}
with the property that
\begin{enumerate}[label=\roman*)]
    \item Each $V[[t]]_{\lambda}$ is a finite rank free $R[[t]]$-module invariant under $\nabla_{t^2\frac{\partial}{\partial t}}$ and
    \item ${V[[t]]_{\lambda}}|_{t=0}=V_{\lambda}$, where $V_{\lambda}$ is the $\lambda$-generalized eigenspace of $A_0$.
\end{enumerate}
We now study the uniqueness of the decomposition (2.9) with the above properties. The key step is the following lemma.
\begin{lemma}
Let $(V[[t]],\nabla_{t^2\frac{\partial}{\partial t}}=t^2\frac{\partial}{\partial t}+A)$ and $(V'[[t]], \nabla'_{t^2\frac{\partial}{\partial t}}=t^2\frac{\partial}{\partial t}+A')$, where $A=\sum_{i\geq 0}A_it^i, A'=\sum_{i\geq 0}A_i't^i\in \mathrm{End}(V)[[t]]$, be two formal connections with quadratic poles, such that $A_0, A_0'$ each has a single eigenvalue $\lambda,\lambda'$, respectively. If $\lambda-\lambda'$ is invertible in $R$, then any morphism of connections $f=\sum_{i\geq 0}f_it^i: (V[[t]],\nabla_{t^2\frac{\partial}{\partial t}})\rightarrow (V'[[t]], \nabla'_{t^2\frac{\partial}{\partial t}})$ is zero.   
\end{lemma}
\emph{Proof}. The condition that $f$ intertwines the connections is equivalent to requiring that 
\begin{equation}
t^2\frac{\partial}{\partial t}f=fA-A'f.   
\end{equation}
This is equivalent to requiring that
\begin{equation}
f_{r}A_0-A'_0f_{r}=\textrm{a linear expression involving}\;f_1,\cdots,f_{r-1}.    
\end{equation}
The statement is then proved by induction on $r$, where we inductively assume the right hand side of (2.11) is zero (which is automatically true for the base case $r=0$), which implies $f_rA_0-A_0'f_r=0$. Rewrite this as 
\begin{equation}
(f_r(A_0-\lambda)-(A_0'-\lambda')f_r)+(\lambda-\lambda')f_r=0.    
\end{equation}
Since (2.7) is an isomorphism, we conclude that $f_r=0$. \qed\par\indent
The following is an immediate corollary of Lemma 2.2.
\begin{cor}
Suppose $(V[[t]],\nabla_{t^2\frac{\partial}{\partial t}}=t^2\frac{\partial}{\partial t}+A)$ and $(V'[[t]], \nabla'_{t^2\frac{\partial}{\partial t}}=t^2\frac{\partial}{\partial t}+A')$ are two formal connections with quadratic poles, and that there are decompositions 
\begin{equation}V[[t]]=\bigoplus_{\lambda} V[[t]]_{\lambda},\quad V'[[t]]=\bigoplus_{\lambda'} V'[[t]]_{\lambda'} \end{equation}
satisfying conditions i) and ii). Assuming that each difference $\lambda-\lambda'$, where $\lambda\in \mathrm{spec}(A_0), \lambda'\in\mathrm{spec}(A_0')$, is invertible in $R$, then any morphism $f: V[[t]]\rightarrow V'[[t]]$ of connections satisfies $f(V[[t]]_{\lambda})\subset V'[[t]]_{\lambda}$ for all $\lambda\in R$ (if $\lambda$ is not an eigenvalue of $A_0$, it is understood that $V[[t]]_{\lambda}=0$).\qed
\end{cor}
\begin{cor}
Assume that the difference between any two distinct eigenvalues of $A_0$ is invertible in $R$. Then for any formal connection $(V[[t]],\nabla_{t^2\frac{\partial}{\partial t}})$ as in (2.1) with quadratic pole term $A_0$, there exists a unique decomposition of the form (2.9) satisfying i) and ii).     
\end{cor}
\emph{Proof}. Existence follows from Lemma 2.1, and uniqueness follows from applying Corollary 2.3 to $f=\mathrm{id}$. \qed\par\indent
We call the decomposition as in Corollary 2.4 \emph{the elementary splitting} of $(V[[t]],\nabla_{t^2\frac{\partial}{\partial t}})$. The following proposition gives a useful characterization of the elementary splitting. 
\begin{prop}
In the setting of Corollary 2.4, let $\{F_{\lambda}\}_{\lambda\in \mathrm{spec}(A_0)}, F_{\lambda}\in \mathrm{End}(V)[[t]]$ be a collection of endomorphisms of the connection $(V[[t]],\nabla_{t^2\frac{\partial}{\partial t}})$ labeled by the eigenvalues of $A_0$. Suppose that
\begin{enumerate}[label=\arabic*)]
    \item $F_{\lambda}^2=F_\lambda$ for all $\lambda$;
    \item $F_{\lambda}\circ F_{\lambda'}=0$ for $\lambda\neq \lambda'$;
    \item $\sum_{\lambda\in \mathrm{spec}(A_0)} F_{\lambda}=\mathrm{id}$;
    \item ${F_{\lambda}}|_{t=0}\in\mathrm{End}(V)$ is the projection onto $V_{\lambda}$.
\end{enumerate}
Then $F_{\lambda}$ is the projection onto $V[[t]]_{\lambda}$ as in the elementary splitting. 
\end{prop}
\emph{Proof}. Assumptions 1)-3) imply that there is a direct sum decomposition 
\begin{equation}
V[[t]]=\bigoplus_{\lambda} \mathrm{im}(F_{\lambda}),
\end{equation}
and $F_{\lambda}$ is the projection onto $\mathrm{im}(F_{\lambda})$. We show that the decomposition (2.14) satisfies the conditions of Corollary 2.4, whence by uniqueness implies that it coincides with the elementary splitting. 
\begin{itemize}
    \item 4) immediately implies condition ii) of Corollary 2.4.
    \item Since $F_{\lambda}$ is an endomorphism of $\nabla_{\frac{\partial}{\partial t}}$, its image is preserved by $\nabla_{\frac{\partial}{\partial t}}$. This shows half of i).
    \item To show the other half of i), we need to prove that $\mathrm{im}(F_{\lambda})$ is a finite rank free $R[[t]]$-module for each $\lambda$. For each $\lambda$, fix an $R$-basis $e^{\lambda}_1,\cdots,e^{\lambda}_{k_{\lambda}}$ of $V_{\lambda}$. We claim that $\{F_{\lambda}(e^{\lambda}_i)\}_{i=1}^{k_{\lambda}}$ form a basis of $\mathrm{im}(F_{\lambda})$. To see this, first note that by 4), ${F_{\lambda}(e^{\lambda}_i)}|_{t=0}=e^{\lambda}_i$. Since $\{e^{\lambda}_i\}$ form a basis of $V^{\lambda}$, there is an injection
    \begin{equation}
    (F_{\lambda}(e^{\lambda}_1),\cdots,F_{\lambda}(e^{\lambda}_{k_{\lambda}})): R[[t]]^{\oplus k_{\lambda}}\hookrightarrow \mathrm{im}(F_{\lambda}).   
    \end{equation}
    Summing over $\lambda$ we get an injection
    \begin{equation}
      \bigoplus_{\lambda}(F_{\lambda}(e^{\lambda}_1),\cdots,F_{\lambda}(e^{\lambda}_{k_{\lambda}})): \bigoplus_{\lambda}R[[t]]^{\oplus k_{\lambda}}\hookrightarrow \bigoplus_{\lambda} \mathrm{im}(F_{\lambda})=V[[t]]. 
    \end{equation}
    However, the $t=0$ part of the collection $\{F_{\lambda}(e^{\lambda}_1),\cdots,F_{\lambda}(e^{\lambda}_{k_{\lambda}})\}_{\lambda\in \mathrm{spec}(A_0)}$ exactly coincides with a generalized eigenbasis of $V$, and thus the map (2.16) is invertible, which further implies each of its individual pieces as in (2.15) is invertible. \qed
\end{itemize} 

\renewcommand{\theequation}{3.\arabic{equation}}
\setcounter{equation}{0}

\section{Quantum Steenrod operations}
In this section, we review the definition of quantum Steenrod operations following \cite{SW} using the Morse chain model of quantum cohomology. We then study its properties relative to the quantum connection and elementary splitting of Section 2. \par\indent
\textbf{Quantum Cohomology}. Fix a ground ring $R$, a Morse function $f$ and a metric $g$ on $X$ such that the associated gradient flow is Morse-Smale. Consider the $\mathbb{Z}$-graded Morse complex $CM^*(f,R[q])$ generated by critical points of $f$ over $R[q]$, where $q$ is a formal variable of degree $2$, with the grading given by the Morse index together with the $q$-grading. The differential is given by counting gradient flow lines between two critical points whose indices differ by $1$. We use the cohomological convention, i.e. $|x|=\dim W^s(x), 2n-|x|=\dim W^u(x)$, where $W^s,W^u$ are the stable and unstable manifold of $x$, respectively. \par\indent
The chain level quantum cup product can be defined as follows, cf. \cite[section 3]{SW}. Given a homology class $A\in H_2(X,\mathbb{Z})$, inputs $x_0,x_1\in\mathrm{crit}(f)$ and output $x_{\infty}\in \mathrm{crit}(f)$,  let $\mathcal{M}_A(C,x_0,x_1,x_{\infty})$ be the moduli space of $J$-holomorphic maps $u:C=\mathbb{C}P^1\rightarrow X$ with three marked points $z_0,z_1,z_{\infty}\in C$ such that $u$ lies in class $A$, and $z_0,z_1,z_{\infty}$ are constrained on $W^u(x_0), W^u(x_1)$ and $W^s(x_{\infty})$, respectively. Equivalently, it is the moduli space of $J$-holomorphic $u: \mathbb{C}P^1\rightarrow X$ together with gradient half-flowlines \begin{equation}
y_0,y_1:(-\infty,0]\rightarrow X,\quad y_{\infty}:[0,\infty)\rightarrow X    
\end{equation}
satisfying
$$y_k'=\nabla f(y_k), y_k(0)=u(z_k), \lim_{s\rightarrow -\infty}y_k(s)=x_k,$$
\begin{equation}
y_{\infty}'=\nabla f(y_{\infty}), y_{\infty}(0)=u(z_{\infty}), \lim_{s\rightarrow \infty}y_{\infty}(s)=x_{\infty}.
\end{equation}
For a generic choice of $J$, the moduli space $\mathcal{M}_A(C,x_0,x_1,x_{\infty})$ is regular of dimension $2c_1(A)+|x_{\infty}|-|x_0|-|x_1|$.\par\indent
We define the \emph{chain level quantum cup product} by the formula
\begin{equation}
 x_0\star_q x_1:=\sum_{\substack{x_{\infty},A:\\2c_1(A)+|x_{\infty}|-|x_0|-|x_1|=0}} |\mathcal{M}_A(C,x_0,x_1,x_{\infty})| q^{c_1(A)}x_{\infty}, 
\end{equation}
where $|\mathcal{M}_A(C,x_0,x_1,x_{\infty})|$ denotes the signed count of isolated elements of $\mathcal{M}_A(C,x_0,x_1,x_{\infty})$. Note that for a fixed $x_{\infty}$, its coefficient in (3.3) is indeed a well-defined element of $R[q]$ because of monotonicity. Extending (3.3) $q$-linearly, we obtain a chain level product $\star_q:CM^*(f,R[q])\otimes CM^*(f,R[q])\rightarrow CM^*(f,R[q])$ which descends to cohomology. This will be called the \emph{quantum cup product} of $X$, and we denote $QH^*(X,R[q])$ for the resulting cohomological algebra $(H^*(X;R)[q],\star_q)$. \par\indent
In this paper, we also consider the restriction to $q=1$ (and hence collapsing the $\mathbb{Z}$-grading to a $\mathbb{Z}/2$-grading) of (3.3), which by an abuse of terminology is also called the quantum cup product, and denoted 
\begin{equation}
\star\;\;(\mathrm{or}\;\star_{q=1}): H^*(X;R)\otimes H^*(X;R)\rightarrow H^*(X;R).    
\end{equation}
\textbf{Quantum Steenrod operations}. Fix a field $\mathbf{k}$ of odd characteristic $p$. To define the quantum Steenrod operations, following \cite[section 4a]{SW}, we consider a moduli problem with fixed domain but parametrized (equivariant) perturbation data. The relevant perturbation data will be parametrized by
\begin{equation}
S^{\infty}:=\{w=(w_0,w_1,\cdots)\in\mathbb{C}^{\infty}: w_k=0\;\textrm{for}\;k\gg 0, \|w\|^2=1\}.
\end{equation}
For a prime $p$, there is a $\mathbb{Z}/p$-action on $S^{\infty}$ where the standard generator $\tau\in\mathbb{Z}/p$ acts by
\begin{equation}
\tau(w_0,w_1,\cdots)=(e^{2\pi i/p} w_0,e^{2\pi i/p}w_1,\cdots).    
\end{equation}
Consider the cells
\begin{equation}
\Delta_{2k}=\{w\in S^{\infty}: w_k\geq 0, w_{k+1}=w_{k+2}=\cdots=0\},
\end{equation}
\begin{equation}
\Delta_{2k+1}=\{w\in S^{\infty}: e^{i\theta}w_k\geq 0\;\textrm{for some}\;\theta\in[0,2\pi/p], w_{k+1}=w_{k+2}=\cdots=0\}.
\end{equation}
We identify the tangent space of $\Delta_{2k}$ at the point $w_{k}=1$ (and hence all other coordinates are zero) with $\mathbb{C}^k$, via the projection onto the first $k$ coordinates; we use the induced complex orientation on $\Delta_{2k}$. The tangent space of $\Delta_{2k+1}$ at the same point is canonically identified with $\mathbb{C}^k\times i\mathbb{R}$, and we use the complex orientation on the first factor and the positive vertical orientation on the second factor. With these chosen orientations, one has
\begin{equation}
\partial\Delta_{2k}=\Delta_{2k-1}+\tau\Delta_{2k-1}+\cdots+\tau^{p-1}\Delta_{2k-1},
\end{equation}
\begin{equation}
\partial\Delta_{2k+1}=\tau\Delta_{2k}-\Delta_{2k}.
\end{equation}
Let $C=\mathbb{C}P^1$, equipped with the $\mathbb{Z}/p$ action given by rotating by $e^{2\pi i/p}$. Denote the action of the generator by $\sigma$. Fix $p+2$ points $z_0=0, z_k=e^{2\pi ik/p}, k=1,2,\cdots,p, z_{\infty}=\infty$ on $C$. \par\indent
Choose perturbation data $\nu^{eq}_w\in \mathrm{Hom}(TC, TX)$ parametrized by $w\in S^{\infty}$ satisfying the equivariance condition
\begin{equation}
\sigma^*\nu^{eq}_w=\nu^{eq}_{\tau(w)}.
\end{equation}
Let $\mathcal{M}_A(\Delta_i\times C,x_0,x_1,\cdots,x_p,x_{\infty})$ denote the moduli space of pairs $(w,u)$, where $w\in \Delta_i$ and $u: C\rightarrow X$ in class $A$ satisfying Floer's equation
\begin{equation}
(du-\nu^{eq}_w)^{0,1}_{J}=0
\end{equation}
and incidence conditions to $W^u(x_0),\cdots,W^u(x_p),W^s(x_{\infty})$. This is a moduli space with virtual dimension
\begin{equation}
\dim\mathcal{M}_A(\Delta_i\times C,x_0,x_1,\cdots,x_p,x_{\infty})=i+2c_1(A)+|x_{\infty}|-|x_0|-|x_1|-\cdots-|x_{p}|.
\end{equation}
For a generic $J$ and equivariant perturbation data $\nu^{eq}$, a standard transversality argument implies that $\mathcal{M}_A(\Delta_i\times C,x_0,x_1,\cdots,x_p,x_{\infty})$ is regular in dimension $0$ and $1$, which is all we will need. 
Moreover, because of the condition imposed by (3.11), we have an identification
\begin{equation}
\mathcal{M}_A(\tau^j(\Delta_i)\times C,x_0,x_1,\cdots,x_p,x_{\infty})\cong \mathcal{M}_A(\Delta_i\times C,x_0,x_{p-j+1},\cdots,x_p,x_1,\cdots,x_{p-j},x_{\infty})
\end{equation}
given by $(w,u)\mapsto (\tau^{-j}(w),u\circ \sigma^{-j})$. This defines maps, for $i\geq 0$,
\begin{equation}
\Sigma_A^{i}:=\Sigma_A(\Delta_i,\cdots):CM^*(f)\otimes CM^*(f)^{\otimes p}\rightarrow CM^{*-i-2c_1(A)}(f).
\end{equation}
Let $t,\theta$ be formal variables with $\mathbb{Z}$-grading $|t|=2, |\theta|=1, t\theta=\theta t,\theta^2=0$. Then $CM(f)[[t,\theta]]$ can be viewed as the $\mathbb{Z}/p$-equivariant complex of $CM^*(f)$, where $\mathbb{Z}/p$ acts trivially on $CM^*(f)$. Fixing a Morse cocycle $b\in CM^*(f)$, one can combine the maps in (3.15) into a chain map
\begin{equation}
\Sigma_{A,b}: CM^*(f)\rightarrow (CM^*(f)[[t,\theta]])^{*+p|b|-2c_1(A)}
\end{equation}
given by
\begin{equation}
x\mapsto (-1)^{|b||x|}\sum_k\big(\Sigma_A(\Delta_{2k},x,b,\cdots,b)+(-1)^{|b|+|x|}\Sigma_A(\Delta_{2k+1},x,b,\cdots,b)\theta\Big)t^k.
\end{equation}
Up to homotopy, (3.17) only depends on the cohomology class $[b]\in QH^*(X;\mathbf{k})$, cf. \cite[Lemma 4.4.]{SW}. Finally, summing over $A$ (which is again well defined by monotonicity), we obtain a chain map
\begin{equation}
\Sigma^q_b:= \sum_{A}q^{c_1(A)}\Sigma_{A,b}:CM^*(f)[q]\rightarrow   (CM^*(f)[[q,t,\theta]])^{*+p|b|}.  
\end{equation}
Finally, extending $(t,\theta)$-linearly, we obtain a chain map 
\begin{equation}
\Sigma^q_b: CM^*(f)[[q,t,\theta]]\rightarrow CM^*(f)[[q,t,\theta]]    
\end{equation} 
of degree $p|b|$. We denote the cohomology level map of $\Sigma^q_b$ as $Q\Sigma^q_b$. As before, we denote by the restriction to $q=1$ of $Q\Sigma^q_b$ as $Q\Sigma$. 
\begin{rmk}
As per our convention, $[[-]]$ always involves taking graded completion. In particular, in the $\mathbb{Z}$-graded setting, because $H^*(X;\mathbf{k})$ is in bounded degrees, in fact $H^*(X;\mathbf{k})[[q,t,\theta]]=H^*(X;\mathbf{k})[q,t,\theta]$. However, when we set $q=1$ (and hence collapses the $\mathbb{Z}$-grading to a $\mathbb{Z}/2$-grading), elements of $H^*(X;\mathbf{k})[[t,\theta]]$ are genuine power series in $(t,\theta)$. Therefore, `restriction to $q=1$' of a certain operation should be understood as `setting $q=1$ and then taking $t$-completion'.     
\end{rmk}

We collect some well-known properties of quantum Steenrod operations, cf. for instance \cite{SW}. First, $Q\Sigma^q$ defines a \emph{Frobenius $p$-linear algebra action} of $QH^*(X,\mathbf{k}[q])$ on $H^*(X;\mathbf{k})[[q,t,\theta]]$. More precisely, 
\begin{itemize}
     \item (Unitality) $Q\Sigma^q_e=\mathrm{id}$.
    \item (Additivity) $Q\Sigma^q_{b+b'}=Q\Sigma^q_b+Q\Sigma^q_{b'}$.
    \item (Frobenius linearity) $Q\Sigma^q_{w b}=w^pQ\Sigma^q_b$ where $w$ is a constant.
    \item (Quantum Cartan relation) 
\begin{equation} 
Q\Sigma^q_b\circ Q\Sigma^q_{b'}=(-1)^{\frac{p(p-1)}{2}|b||b'|} Q\Sigma^q_{b\star_q b'}.
\end{equation}
\end{itemize}
By restricting to $q=1$, $Q\Sigma$ also satisfies analogous properties of (3.20). Furthermore, the classical ($q=0$) and non-equivariant ($t=\theta=0$) parts of $Q\Sigma^q_b$ are respectively
\begin{itemize}
    \item ${Q\Sigma^q_b|}_{q=0}=St(b)\cup-$, where $St$ denotes the classical Steenrod powers.
    \item ${Q\Sigma^q_{b}}|_{t=\theta=0}=b^{\star p}\star_q-$.
\end{itemize}
For the purpose of this paper, it suffices to consider quantum Steenrod operations restricted to the $S^1$-equivariant quantum cohomology. More precisely, fixing $b\in QH^*(X;\mathbf{k}[q])$, we consider the composition
\begin{equation}
 H^*(X;\mathbf{k})[[q,t]]\hookrightarrow H^*(X;\mathbf{k})[[q,t,\theta]]\xrightarrow{Q\Sigma_b} H^*(X;\mathbf{k})[[q,t,\theta]]\xrightarrow{\theta=0} H^*(X;\mathbf{k})[[q,t]].  
\end{equation}
By an abuse of notation, we also denote (3.21) as $Q\Sigma_b^q$. It inherits the prior-mentioned properties of the usual quantum Steenrod operations and defines a Frobenius $p$-linear algebra action of $QH^*(X;\mathbf{k}[q])$ on $H^*(X;\mathbf{k})[[q,t]]$.\par\indent 
We now recall the main result of \cite{SW}, called the \emph{covariant constancy property of quantum Steenrod operations}. Before stating the theorem, we first recall the definition of the \emph{quantum $q$-connection}, which is defined by (it is convenient to view this connection as differentiating in the $tq\frac{\partial}{\partial q}$-direction)
\begin{equation}
 \nabla^{QH}_{tq\frac{\partial}{\partial q}}:=tq\frac{\partial}{\partial q}+c_1\star_q: H^*(X;R)[[q,t]]\rightarrow  H^*(X;R)[[q,t]].
\end{equation}
\begin{thm} (\cite[Theorem 1.4]{SW}) For any $b\in QH^*(X;\mathbf{k}[q])$,
\begin{equation}
Q\Sigma^q_b \circ \nabla^{QH}_{tq\frac{\partial}{\partial q}}- \nabla^{QH}_{tq\frac{\partial}{\partial q}}\circ Q\Sigma^q_b=0.
\end{equation}   
\end{thm}
The following proposition relates the $p$-curvature of the quantum $t$-connection with quantum Steenrod operations. The idea that these two objects are related originated from the work of Jae Hee Lee \cite{Lee}, which studied the $p$-curvature of quantum connections of symplectic resolutions arising from representation theory. The adaptation of this idea to the closed monotone setting is due to Paul Seidel, and we thank him for communicating the following result. 
\begin{prop}[Seidel]
Let $X$ be a closed monotone symplectic manifold. Then over any field $\mathbf{k}$ of odd characteristic $p$,  
\begin{equation}
Q\Sigma_{c_1}+F^{QH}_{t^2\frac{\partial}{\partial t}}
\end{equation}
is a nilpotent operation on $H^*(X;\mathbf{k})[[t]]$.
\end{prop}
\emph{Proof}. We first show that
\begin{equation}
Q\Sigma^q_{c_1}-F^{QH}_{tq\frac{\partial}{\partial q}}
\end{equation}
is nilpotent, where we recall that
\begin{equation}
F^{QH}_{tq\frac{\partial}{\partial q}}:=(\nabla^{QH}_{tq\frac{\partial}{\partial q}})^p-t^{p-1}\nabla^{QH}_{tq\frac{\partial}{\partial q}}
\end{equation}
is the $p$-curvature of the quantum $q$-connection along the vector field $tq\frac{\partial}{\partial q}$, cf. (A.1). The proof is essentially a consequence of Theorem 3.2. \par\indent
First we note that the classical ($q=0$) terms of both $Q\Sigma^q_{c_1}$ and $F^{QH}_{tq\frac{\partial}{\partial q}}$ are the classical Steenrod operation $(-t^{p-1}c_1+c_1^{\cup p})\cup$. Moreover, both commute with $\nabla^{QH}_{tq\frac{\partial}{\partial q}}$: for $Q\Sigma^q_{c_1}$ this is Theorem 3.2 and for $F^{QH}_{tq\frac{\partial}{\partial q}}$ it is obvious as the $p$-curvature is just a combination of the connection. This implies that their difference is $O(q^p)$. To see this, we write 
\begin{equation}
Q\Sigma^q_{c_1}-F^{QH}_{tq\frac{\partial}{\partial q}}=A_1(t)q+A_2(t)q^2+\cdots,   
\end{equation}
where $A_i$'s are polynomial expressions in $t$ with coefficients in $\mathrm{End}(H^*(X;\mathbf{k}))$. Rewriting 
\begin{equation}
\nabla^{QH}_{tq\frac{\partial}{\partial q}}=tq\frac{\partial}{\partial q}+(c_1\cup+O(q))
\end{equation}
and equating powers of $q$ in the expression $[Q\Sigma^q_{c_1}-F^{QH}_{tq\frac{\partial}{\partial q}},\nabla^{QH}_{tq\frac{\partial}{\partial q}}]=0$, one obtains that for all $r$,
\begin{equation}
rtA_r+[A_r,c_1\cup]=\textrm{linear expression in}\;\;\{A_1,\cdots,A_{r-1}\}.   
\end{equation}
We can argue by induction (in the range $0\leq r<p$) and assume the right hand side of (3.29) is zero. Then, since $r<p$, it is invertible in $\mathbf{k}$. On the other hand, $[-,c_1\cup]$ is nilpotent. This implies that $A_r=0$.\par\indent 
Moreover, the non-equivariant ($t=0$) terms of both $Q\Sigma^q_{c_1}$ and $F^{QH}_{tq\frac{\partial}{\partial q}}$ are equal to $c_1^{\star p}\star$. So their difference is in fact of the form
\begin{equation}
(Q\Sigma^q_{c_1}-F^{QH}_{tq\frac{\partial}{\partial q}})(x)=y_1tq^p+\textrm{terms with higher powers of $(t,q)$}.
\end{equation}
Since the total degrees of both $Q\Sigma^q_{c_1}$ and $F^{QH}_{tq\frac{\partial}{\partial q}}$ are $2p$, and the degrees of $(t,q)$ are $(2,2)$, respectively, the coefficients $y_i$'s have degree at most $\deg(x)-2$. Since $H^*(X;\mathbf{k})$ lies in bounded degree, we conclude that $Q\Sigma^q_{c_1}-F^{QH}_{tq\frac{\partial}{\partial q}}$ is nilpotent.  \par\indent
To conclude the proof, we note that by Lemma B.2
\begin{equation}
F^{QH}_{tq\frac{\partial}{\partial q}}=-F^{QH}_{t^2\frac{\partial}{\partial t}}\quad\mathrm{on}\;\;H^*(X;\mathbf{k})[[q,t]].
\end{equation}
Restricting to $q=1$ gives the desired result.\qed

\renewcommand{\theequation}{4.\arabic{equation}}
\setcounter{equation}{0}
\section{Proof of Theorem 1.3}
We first reduce the proof of Theorem 1.3 to a computation of the $p$-curvature of the quantum connection, using Katz's local monodromy theorem, cf. Theorem A.1. Fix a subring $R$ of $\overline{\mathbb{Q}}$ which is finitely generated (as an algebra) over $\mathbb{Z}$, with the additional properties that
\begin{enumerate}[label=\arabic*)]
    \item $c_1\star \in\mathrm{End}(H^*(X;R))$ admits a Jordan decomposition, cf. beginning of Section 2.
    \item The difference between two distinct eigenvalues of $c_1\star$ is invertible in $R$. 
    \item $2$ is invertible in $R$ (the grading operator $\mu$ in the quantum $t$-connection, cf. Definition 1.1, involves $\frac{1}{2}$). 
\end{enumerate}
Such a choice exists for the following elementary reasons.
\begin{itemize}
    \item By the universal coefficient theorem and compactness of $X$, there exists a localization of $\mathbb{Z}$ at finitely many elements, denoted $R$, such that $H^*(X;R)$ is a finite rank free $R$-module. 
    \item By monotonicity, the Gromov-Witten invariants of $X$ are defined over $\mathbb{Z}$ (and hence $R$). Adjoin to $R$ the finitely many eigenvalues $\{\lambda_1,\cdots,\lambda_m\}$ of $c_1\star\in\mathrm{End}(H^*(X;R))$ in some algebraic closure of $\mathrm{Frac}(R)$, and invert their pairwise (nonzero) differences. Call the result $R'$. 
    \item There are polynomials $(p_2(t),\cdots,p_m(t); q_2(t),\cdots,q_m(t))$ with coefficients in $\mathrm{Frac}(R')$ such that $1=p_i(t)(t-\lambda_i)^{k_i}+q_i(t)\prod_{j=1}^{i-1}(t-\lambda_j)^{k_j}$ for $i=2,\cdots,m$. Adjoining all the coefficients of $\{p_i(t), q_i(t)\}$ to $R'$, we obtain a new ring $R''$ over which there is a direct sum decomposition 
    \begin{equation}
      H^*(X;R'')=\bigoplus_{i=1}^m \mathrm{ker}(c_1\star-\lambda_iI)^{k_i}.  
    \end{equation}
    \item Finally, since $R''$ is finitely generated over $\mathbb{Z}$, it is Noetherian. Hence each $\mathrm{ker}(c_1\star-\lambda_iI)^{k_i}$ is finitely presented. Since local freeness of a finitely presented module is a stalk local property, there exists $f_i\in R''$ such that $\mathrm{ker}(c_1\star-\lambda_iI)^{k_i}_{f_i}$ is free. In particular, after base changing to $R''':=R''_{f_1\cdots f_m}$, each summand of the decomposition (4.1) becomes a finite rank free $R'''$-module. Replace $R$ by $R'''$ gives the desired ring.
\end{itemize}
Since the quadratic pole term of (1.2) is $-c_1\star$, to keep the notation consistent with Section 2, we denote the $\lambda$-generalized eigenspace of $c_1\star$ by $H^*(X;R)_{-\lambda}$. Recall from Corollary 2.4 that there is a unique finite decomposition of $\nabla^{QH}_{\frac{\partial}{\partial t}}$ into finite rank free $R[[t]]$-modules with connections
\begin{equation}
\nabla^{QH}_{\frac{\partial}{\partial t}}=\bigoplus_{\lambda\in\mathrm{spec}(c_1\star)} \nabla^{QH,-\lambda}_{\frac{\partial}{\partial t}},    
\end{equation}
where up to a gauge transformation, 
\begin{equation}
\nabla^{QH,-\lambda}_{\frac{\partial}{\partial t}}=\frac{\partial}{\partial t}-\frac{\lambda I+N_{\lambda}}{t^2}+\sum_{m\geq -1}A_{m} t^m
\end{equation}
with $N_{\lambda}$ a nilpotent matrix. We can rewrite this as
\begin{equation}
\nabla^{QH,-\lambda}_{\frac{\partial}{\partial t}}=\mathcal{E}^{\frac{-\lambda}{t^2}}\otimes\tilde{\nabla}^{QH,-\lambda}_{\frac{\partial}{\partial t}}
\end{equation}
where $\mathcal{E}^{-\frac{\lambda}{t^2}}=(R((t)),\frac{\partial}{\partial t}-\frac{\lambda}{t^2})$ and
\begin{equation}
\tilde{\nabla}^{QH,-\lambda}_{\frac{\partial}{\partial t}}=\frac{\partial}{\partial t}-\frac{N_{\lambda}}{t^2}+\sum_{m\geq -1}A_{m} t^m,
\end{equation}
which we call the \emph{residual connection}. 
%\begin{mydef}
%The connection $\nabla^{QH}_{\frac{\partial}{\partial t}}$ is said to be of \emph{unramified exponential type} if for each $\lambda$, $\tilde{\nabla}^{QH}_{\lambda}$ of (6.4) has a regular singularity at $t=0$.  %\end{mydef}
By definition of exponential type, Theorem 1.3 is an immediate consequence of the following proposition.
\begin{prop}
For each $\lambda\in \mathrm{spec}(c_1\star)$, the residual connection $(H^*(X;R)((t))_{-\lambda},\tilde{\nabla}^{QH,-\lambda}_{\frac{\partial}{\partial t}})$ of (4.5) has regular singularity at $t=0$ and quasi-unipotent monodromy after based-changed to $K((t))$, where $K=\mathrm{Frac}(R)$. 
\end{prop}
By Theorem A.1, Proposition 4.1 would follow if for all $\mathfrak{m}\in \mathrm{mSpec}(R)$, the $p$-curvature of $\tilde{\nabla}^{QH}_{\lambda}\otimes_{R((t))}\kappa(\mathfrak{m})((t))$ is nilpotent (where $\kappa(\mathfrak{m}):=R/\mathfrak{m}$ is the residue field), or equivalently if the following proposition holds, since the $p$-curvature of $\nabla^{QH,-\lambda}_{t^2\frac{\partial}{\partial t}}$ and $\tilde{\nabla}^{QH,-\lambda}_{t^2\frac{\partial}{\partial t}}$ differ by $-\lambda^p$.
\begin{prop}
For each $\lambda\in \mathrm{spec}(c_1\star)$ and $\mathfrak{m}\in\mathrm{mSpec}(R)$, 
\begin{equation}
F^{{-\lambda}}_{t^2\frac{\partial}{\partial t}}+\lambda^p    
\end{equation}
is nilpotent, where $F^{{-\lambda}}_{t^2\frac{\partial}{\partial t}}$ is the $p$-curvature of $\nabla^{QH}_{\lambda}\otimes_{R((t))}\kappa(\mathfrak{m})((t))$ along the vector field $t^2\frac{\partial}{\partial t}$:
\begin{equation}
F^{{-\lambda}}_{t^2\frac{\partial}{\partial t}}:=(\nabla^{QH,-\lambda}_{t^2\frac{\partial}{\partial t}})^p=(t^2\frac{\partial}{\partial t}-(\lambda I+N_{\lambda})+\sum_{m\geq -1}A_{m} t^{m+2})^p.  
\end{equation} 
\end{prop}
The proof of Proposition 4.2 requires two ingredients. The first of which is Proposition 3.3. The second is the following compatibility between $Q\Sigma$ and the elementary splitting of $\nabla^{QH}_{\frac{\partial}{\partial t}}$, cf. Section 2.\par\indent
Let $\mathbf{k}$ be a field of odd characteristic $p$ satisfying conditions 1) - 3) at the beginning of this section. Let $H^*(X;\mathbf{k})[[t]]_{-\lambda}$ denote the $(-\lambda)$-component of the unique elementary splitting (cf. Corollary 2.4). 
\begin{lemma}
Let $e_{\lambda}\in H^*(X;\mathbf{k})_{-\lambda}$ be the $\lambda$-idempotent, i.e. the projection of the unit $e$ onto the $\lambda$-generalized eigenspace of $c_1\star$. Then, $Q\Sigma_{e_{\lambda}}$ is equal to the projection onto $H^*(X;\mathbf{k})[[t]]_{-\lambda}$.
\end{lemma}
\emph{Proof}. We consider the collection $\{Q\Sigma_{e_{\lambda}}\}_{\lambda}$ and show that they satisfy the assumptions of Proposition 2.5. First, $Q\Sigma_{e_{\lambda}}$ is covariantly constant with respect to $\nabla^{QH}_{\frac{\partial}{\partial t}}$ by Lemma B.2. Next we check conditions 1)-4) of Proposition 2.5. By the quantum Cartan relations, 
\begin{equation}
Q\Sigma_{e_{\lambda}}^2=Q\Sigma_{e_{\lambda}\star e_{\lambda}}= Q\Sigma_{e_{\lambda}}   
\end{equation}
and
\begin{equation}
Q\Sigma_{e_{\lambda}}\circ Q\Sigma_{e_{\lambda'}}=Q\Sigma_{e_{\lambda}\star e_{\lambda'}}=0\quad\mathrm{if}\;\;\lambda\neq \lambda'
\end{equation}
since two distinct idempotents are orthogonal with respect to the quantum cup product. This verifies 1) and 2). 3) follows from the additivity and unitality of quantum Steenrod operations and that $\sum_{\lambda} e_{\lambda}=e$. Finally, to prove 4), recall that the non-equivariant part of quantum Steenrod operations is given by
\begin{equation}
{Q\Sigma_{e_{\lambda}}}|_{t=0}=e_{\lambda}^{\star p}\star=e_{\lambda}\star=\textrm{projection onto}\;\;H^*(X;\mathbf{k})_{-\lambda}.    
\end{equation}
Lemma 4.3 then follows from Proposition 2.5.\qed

\begin{cor}
If $b\in H^*(X;\mathbf{k})_{-\lambda}$ and $y\in H^*(X;\mathbf{k})[[t]]_{-\lambda'}$ with $\lambda\neq \lambda'$, then $Q\Sigma_b(y)=0$.  
\end{cor}
\emph{Proof}. By the quantum Cartan relations,
\begin{equation}
Q\Sigma_b=Q\Sigma_{b\star e_{\lambda}}=Q\Sigma_b\circ Q\Sigma_{e_{\lambda}}.
\end{equation}
By Lemma 4.3, $Q\Sigma_{e_{\lambda}}$ is the projection onto $H^*(X;\mathbf{k})[[t]]_{-\lambda}$, and hence is zero when applied to $y\in H^*(X;\mathbf{k})[[t]]_{-\lambda'}, \lambda\neq \lambda'$.\qed\par\indent
\emph{Proof of Proposition 4.2}.  Fix an arbitrary $\mathfrak{m}\in\mathrm{mSpec}(R)$. Over $\kappa(\mathfrak{m})$, both $Q\Sigma_{c_1}$ and $F^{QH}_{t^2\frac{\partial}{\partial t}}$ respects the elementary splitting by Corollary 2.3. Therefore, Proposition 3.3 implies that for each $\lambda\in\mathrm{spec}(c_1\star\in \mathrm{End}(H^*(X;R)))$ (by abuse of notation, we also denote by $\lambda$ its image in the residue field $\kappa(\mathfrak{m})$, which gives an element of $\mathrm{spec}(c_1\star\in \mathrm{End}(H^*(X;\kappa(\mathfrak{m}))))$), 
\begin{equation}
(Q\Sigma_{c_1}+F^{QH}_{t^2\frac{\partial}{\partial t}})|_{H^*(X;\kappa(\mathfrak{m}))[[t]]_{-\lambda}}=(Q\Sigma_{c_1})|_{H^*(X;\kappa(\mathfrak{m}))[[t]]_{-\lambda}}+F^{-\lambda}_{t^2\frac{\partial}{\partial t}}
\end{equation}
is nilpotent. By Lemma B.2, $Q\Sigma_{c_1}$ commutes with $\nabla^{QH}_{t^2\frac{\partial}{\partial t}}$, and in particular it commutes with $F^{QH}_{t^2\frac{\partial}{\partial t}}$. Therefore, Proposition 4.2 would follow if we can show that 
\begin{equation}
(Q\Sigma_{c_1}-\lambda^p)|_{H^*(X;\kappa(\mathfrak{m}))[[t]]_{-\lambda}}    
\end{equation}
is nilpotent. \par\indent  
To show that (4.13) is nilpotent, we write $c_1=\sum_{\lambda\in\mathrm{spec}(c_1\star)}c_1^{\lambda}$, where $c_1^{\lambda}$ is the projection of $c_1$ onto $H^*(X;\kappa(\mathfrak{m}))_{-\lambda}$. First we note that $c_1^{\lambda}-\lambda e_{\lambda}$ is a nilpotent element of $H^*(X;\kappa(\mathfrak{m}))_{-\lambda}\subset H^*(X;\kappa(\mathfrak{m}))$. Indeed, by definition of a generalized eigenspace, there exists an integer $k$ such that $(c_1-\lambda e)^{\star k}\star=(c_1^{\lambda}-\lambda e_{\lambda})^{\star k}\star =0$ on $H^*(X;\kappa(\mathfrak{m}))_{-\lambda}$, and thus 
\begin{equation}
(c_1^{\lambda}-\lambda e_{\lambda})^{\star k}=(c_1^{\lambda}-\lambda e_{\lambda})^{\star k}\star e=(c_1^{\lambda}-\lambda e_{\lambda})^{\star k}\star e_{\lambda}=0.    
\end{equation}
Since $Q\Sigma$ is a Frobenius-linear algebra action of $QH^*(X;\kappa(\mathfrak{m}))$ on $H^*(X;\kappa(\mathfrak{m}))[[t]]$, we conclude that
\begin{equation}
(Q\Sigma_{c_1^{\lambda}}-\lambda^p Q\Sigma_{e_{\lambda}})^k=Q\Sigma_{(c_1^{\lambda}-\lambda e_{\lambda})^{\star k}}=0.
\end{equation}
By Corollary 4.4,
\begin{equation}
(Q\Sigma_{c_1}-\lambda^pQ\Sigma_e)^k|_{H^*(X;\kappa(\mathfrak{m}))[[t]]_{-\lambda}}=(Q\Sigma_{c_1^{\lambda}}-\lambda^p Q\Sigma_{e_{\lambda}})^k|_{H^*(X;\kappa(\mathfrak{m})[[t]]_{-\lambda}}\stackrel{(4.15)}{=} 0.
\end{equation}
Since $Q\Sigma_e=\mathrm{id}$ by unitality, this shows that (4.13) is nilpotent, which concludes the proof.\qed

\begin{appendices}
\renewcommand{\theequation}{A.\arabic{equation}}
\setcounter{equation}{0}
\section{Katz's local monodromy theorem}
Let $\mathbf{k}$ be a field of characteristic $p$. Let $(M,\nabla)$ be a finite rank free $\mathbf{k}[[t]]$ (or $\mathbf{k}((t))$) module equipped with a connection, and $D\in \mathrm{Der}_{\mathbf{k}}(\mathbf{k}[[t]])$. Recall that the \emph{$p$-curvature} of $\nabla$ along $D$, denoted $F^{\nabla}_{D}$, is defined as
\begin{equation}
F^{\nabla}_D:=\nabla^p_{D}-\nabla_{D^p}.
\end{equation}
In this section, we reproduce a proof of the following version of Katz's local monodromy theorem.
\begin{thm}[Katz's local monodromy theorem (formal version)]
Let $R\subset \overline{\mathbb{Q}}$ be a finitely generated ring extension of $\mathbb{Z}$. %Let $R((t))$ denote the ring of formal formal Laurent series over $R$. 
Let $(M,\nabla)$ be a finite rank free $R((t))$-module equipped with a connection. Suppose that for all maximal ideals $\mathfrak{m}\in\mathrm{mSpec}R$, the connection $(M,\nabla)\otimes_{R((t))}\kappa(\mathfrak{m})((t))$ has nilpotent $p$-curvature, where $p$ is the characteristic of the finite residue field $\kappa(\mathfrak{m})$, then the connection $(M,\nabla)\otimes_{R((t))}\mathrm{Frac}(R)((t))$ has regular singularity and quasi-unipotent monodromy at $t=0$. 
%In fact, after inverting finitely many elements of $R$, $(M,\nabla)$ has regular singularity with quasi-unipotent monodromy at $t=0$. 
\end{thm}
The proof of Theorem A.1 follows almost verbatim from Katz's original proof \cite[Theorem 13.0]{Ka1}, with one minor difference: in \cite{Ka1}, Katz worked over rational functional field $\mathrm{Frac}(R)(t)$ instead of formal Laurent series $\mathrm{Frac}(R)((t))$. The caveat is that, unlike the case of $\mathrm{Frac}(R)(t)$, the denominators in the coefficients of an element of $\mathrm{Frac}(R)((t))$ form an a priori infinite subset of $R$, which can cause potential problems when one tries to reduce to positive characteristics. Fortunately, this can be circumvented by the following lemma, which also follows from a result of Katz. 
\begin{lemma}
In the setting of Theorem A.1, there exists an element $h\in R$ such that $(M,\nabla)\otimes_{R((t))}R[\frac{1}{h}]((t))$ contains a cyclic vector.     
\end{lemma}
Recall that given a finite free $R((t))$-module with integrable connection $(M,\nabla)$, an element $v\in M$ is \emph{cyclic} if $\{v,\nabla_{\frac{\partial}{\partial t}}v,\cdots,\nabla_{\frac{\partial}{\partial t}}^{n-1}v\}$ (or equivalently $\{v,\nabla_{t\frac{\partial}{\partial t}}v,\cdots,\nabla_{t\frac{\partial}{\partial t}}^{n-1}v\}$) form a basis of $M$ over $R((t))$, where $n$ is the rank of $M$. \par\indent
\emph{Proof of Lemma A.2}. \cite[Theorem 2/Remarks (6)]{Ka2} implies that there exist elements $g_0(t),g_1(t)\cdots,g_{n(n-1)}(t)\in R[\frac{1}{(n(n-1))!}]((t))$ (in loc.cit., $g_i(t)$ is denoted $P(t-a_i)$) which generate the unit ideal such that for each $0\leq i\leq n(n-1)$ with $g_i(t)$ nonzero, $(M,\nabla)$ contains a cyclic vector after based changed to $R[\frac{1}{(n(n-1))!}]((t))[\frac{1}{g_i(t)}]$. For such an $i$, since $R[\frac{1}{(n(n-1))!}]((t))[\frac{1}{g_i(t)}]\subset R[\frac{1}{(n(n-1))!}, \frac{1}{h_i}]((t))$, where $h_i$ is the lowest nonzero coefficient of $g_i(t)$, the statement follows by taking $h=(n(n-1))!h_i$. \qed\par\indent
From now on, by base-changing to $R[\frac{1}{h}]((t))$, we will assume without loss of generality that $(M,\nabla)$ contains a cyclic vector. From this point on, the proof of the formal version of the local monodromy theorem proceeds identically to Katz's original proof, cf. \cite[Theorem 11.9]{Ka1},\cite[Theorem 13.0]{Ka1}, which we recall below for the reader's convenience.

\begin{thm}[Fuchs, Turrittin, Lutz]
Let $(M,\nabla)$ be a free $R((t))$-module of rank $n$ with integrable connection that contains a cyclic vector $v\in M$. Suppose $(M,\nabla)\otimes_{R((t))}\mathrm{Frac}(R)((t))$ does not have a regular singularity at $t=0$. Then for every multiple $a$ of $n!$ which is invertible in $R$, after base changing along $t\mapsto t^a$ there exists an $R((t))$-basis $\mathbf{f}$ of $M$ such that
\begin{equation}
\nabla_{t\frac{\partial}{\partial t}}\mathbf{f}=t^{-N}(A+tB)\mathbf{f} 
\end{equation}
for some integer $N\geq 1$, where $A\in \mathbf{M}_n(R)$ is non-nilpotent and $B\in \mathbf{M}_n(R[[t]])$. 
\end{thm}
\emph{Proof of Theorem A.3}. Suppose $\nabla$ does not have regular singularity at $t=0$. Consider the basis $\mathbf{e}:=\{v,\nabla_{t\frac{\partial}{\partial t}} v, \cdots,\nabla^{n-1}_{t\frac{\partial}{\partial t}}v\}$ generated by the cyclic vector $v$. After the base change $t\mapsto t^a$, the connection $\nabla$ can be written in terms of $\mathbf{e}$ as 
\begin{equation}
\frac{1}{a}\nabla_{t\frac{\partial}{\partial t}}\mathbf{e}=\begin{pmatrix}
0&1& &\cdots& &0\\
0&0&1&\cdots&&0\\
\vdots&&&\cdots&&\vdots\\
0&&&\cdots&&1\\
f_0&f_1&&\cdots&&f_{n-1}
\end{pmatrix}\mathbf{e},
\end{equation}
where some $f_i$ has order of vanishing at $t=0$ given by a negative multiple of $n!$. Let $\mathrm{ord}_t(f_i)$ denote the order of vanishing, and consider the positive integer
\begin{equation}
 \nu:=   \max_{0\leq i\leq n-1} \frac{-\mathrm{ord}_t(f_i)}{n-i}.
\end{equation}
Then, under the basis 
\begin{equation}
\mathbf{f}:=   \begin{pmatrix}
1&0&\cdots&&0\\
0&t^{\nu}&\cdots&&0\\
0&0&t^{2\nu}&\cdots&0\\
\vdots&&\cdots&&\vdots\\
0&0&\cdots&&t^{(n-1)\nu}
\end{pmatrix} \mathbf{e},
\end{equation}
the connection matrix of $\frac{1}{a}\nabla_{t\frac{\partial}{\partial t}}$ has the form
\begin{equation}
B=t^{-\nu}\begin{pmatrix}0&0&\cdots&&0\\
0&\nu t^{\nu}&\cdots&&0\\
0&0&2\nu t^{\nu}&\cdots&0\\
\vdots&&\cdots&&\vdots\\
0&0&\cdots&&(n-1)\nu t^{\nu}
\end{pmatrix}  + t^{-\nu}\begin{pmatrix}
0&1& &\cdots& &0\\
0&0&1&\cdots&&0\\
\vdots&&&\cdots&&\vdots\\
0&&&\cdots&&1\\
t^{n\nu}f_0&t^{(n-1)\nu}f_1&&\cdots&& t^{\nu}f_{n-1}
\end{pmatrix}.
\end{equation}
By definition of $\nu$ (A.4), we have $\mathrm{ord}_t(t^{(n-i)\nu}f_i)\geq 0$ for all $0\leq i\leq n-1$ and moreover for at least one value of $i$ equality is achieved. Let $B_{-\nu}:=t^{\nu} B$, then one has 
\begin{equation}
\det (TI_n-B_{-\nu})\equiv\det \begin{pmatrix}
T&-1& &\cdots& &0\\
0&T&-1&\cdots&&0\\
\vdots&&&\cdots&&\vdots\\
0&&&\cdots&&-1\\
-t^{n\nu}f_0&-t^{(n-1)\nu}f_1&&\cdots&&T- t^{\nu}f_{n-1}
\end{pmatrix}\equiv T^n-\sum_{0\leq i\leq n-1} t^{(n-i)\nu}f_iT^i\quad\textrm{modulo}\;\;t.
\end{equation}
Since $\mathrm{ord}_t(t^{(n-i)\nu}f_i)=0$ for some $i$, from (A.7) we obtain that $\det (TI_n-B_{-\nu})$ is not equal to $T^n$ modulo $t$, and in particular $B_{-\nu}$, and thus $B$, cannot be nilpotent. \qed \par\indent
\emph{Proof of Theorem A.1}. Suppose $\nabla$ is not regular singular at $t=0$. Then after replacing $t\mapsto t^a$ (which preserves nilpotence of $p$-curvature), let $\mathbf{f}, A, B$ be specified as in the statement of Theorem A.3.\par\indent
For $i\geq 1$,
\begin{equation}
(\nabla_{t\frac{\partial}{\partial t}})^i\mathbf{f}=t^{-Ni}(A^i+tB_i)\mathbf{f}\quad\textrm{for some}\;\;B_i\in\mathbf{M}_n(R[[t]]).
\end{equation}
By assumption, after passing to the residue field $\kappa(\mathfrak{m})$ of each $\mathfrak{m}\in\mathrm{mSpec}(R)$, there exists some integer $\alpha(\mathfrak{m})$ such that
\begin{equation}
\Big((\nabla_{t\frac{\partial}{\partial t}})^p-\nabla_{t\frac{\partial}{\partial t}}\Big)^{\alpha(\mathfrak{m})}M\subset \mathfrak{m}M,
\end{equation}
where $p$ is the characteristic of $\kappa(\mathfrak{m})$. Plugging (A.8) (with $i=p$) into (A.9) and equating terms with the lowest power of $t$, one deduces that
\begin{equation}
A^{p\cdot\alpha(\mathfrak{m})}\in \mathfrak{m}\mathbf{M}_n(R)\quad\textrm{for all}\;\mathfrak{m}.   
\end{equation}
In particular, the characteristic polynomial $p_A(T)\in R[T]$ of $A$ equals $T^n$ modulo all maximal ideals of $R$, which implies $p_A(T)=T^n$. This implies that $A$ is nilpotent, which is a contradiction. This completes the first part of the theorem. \par\indent
As we have proved that $\nabla$ is regular singular, let $\mathbf{m}$ be a basis of $M$ such that 
\begin{equation}
\nabla_{t\frac{\partial}{\partial t}}\mathbf{m}=(A+tB)\mathbf{m},    
\end{equation}
where $A\in \mathbf{M}_n(R)$ and $B\in\mathbf{M}_n(R[[t]])$. Decompose $A=D+N$, where $D$ is diagonal, $N$ is nilpotent super-triangular and $[D,N]=0$. Passing to the residue field $\kappa(\mathfrak{m})$ of some $\mathfrak{m}\in\mathrm{mSpec}(R)$ and equating the diagonal part of the constant term of (A.9) implies that 
\begin{equation}
(D^p-D)^{\alpha(\mathfrak{m})}\in \mathfrak{m}\mathbf{M}_n(R).    
\end{equation}
In particular, each diagonal entry $d$ of $D$, modulo $\mathfrak{m}$, must lie in $\mathbb{F}_p\subset \kappa(\mathfrak{m})$. Since this holds for all $\mathfrak{m}\in\mathrm{mSpec}(R)$, by the Chebotarev density theorem this implies that $d\in R\cap\mathbb{Q}$, which proves the second part of Theorem A.1. \qed

\renewcommand{\theequation}{B.\arabic{equation}}
\setcounter{equation}{0}
\section{From the $q$-connection to the $t$-connection}
We recall the well-known relation between the quantum $q$-connection (3.22) and the quantum $t$-connection (1.2), and use it to study the relation between their $p$-curvatures.\par\indent
Define the \emph{total degree operator} $\mathrm{Deg}: H^*(X;R)[[q,t]]\rightarrow H^*(X;R)[[q,t]]$ by
\begin{equation}
\mathrm{Deg}:=2(q\frac{\partial}{\partial q}+t\frac{\partial}{\partial t}+\mu).  
\end{equation}
As its name suggests, the effect of $\mathrm{Deg}$ applied to an element $\beta\in H^*(X;R)[[q,t]]$ is to multiply $\beta$ by its total $\mathbb{Z}$-grading (i.e. combining the grading from $H^*(X)$, from $q$ and from $t$), shifted by $n=\dim_{\mathbb{C}}X$. From the explicit formulae, one easily sees that
\begin{equation}
\nabla^{QH}_{t^2\frac{\partial}{\partial t}}=\frac{1}{2}t\,\mathrm{Deg}-\nabla^{QH}_{tq\frac{\partial}{\partial q}},       
\end{equation}
where 
\begin{equation}
\nabla^{QH}_{t^2\frac{\partial}{\partial t}}:=t^2\frac{\partial}{\partial t}+t\mu-c_1\star_q
\end{equation}
denotes the quantum $t$-connection on  $H^*(X;R)[[q,t]]$ (i.e. (1.2) is the restriction of (B.3) to $q=1$, up to a factor of $t^2$). \par\indent
Before stating the next lemma, recall that over a field of characteristic $p$ the following identities for $p$-th powers of vector fields hold: $(\frac{\partial}{\partial t})^p=(t^2\frac{\partial}{\partial t})^p=0, (t\frac{\partial}{\partial t})^p=t\frac{\partial}{\partial t}, (tq\frac{\partial}{\partial q})^p=t^pq\frac{\partial}{\partial q}$. This allows one to simplify the expressions for the relevant $p$-curvatures. 
\begin{lemma}
Fix $R=\mathbf{k}$ a field of odd characteristic $p$. Let $F^{QH}_{t^2\frac{\partial}{\partial t}}$ denote the $p$-curvature of the quantum $t$-connection (B.3) along $t^2\frac{\partial}{\partial t}$, and let $F^{QH}_{tq\frac{\partial}{\partial q}}$ denote the $p$-curvature of the quantum $q$-connection (3.22) along $tq\frac{\partial}{\partial q}$. Then,
\begin{equation}
F^{QH}_{t^2\frac{\partial}{\partial t}}=-F^{QH}_{tq\frac{\partial}{\partial q}}.
\end{equation}
\end{lemma}
\emph{Proof.} This is a straightforward computation. 
\begin{align}
F^{QH}_{t^2\frac{\partial}{\partial t}}&=t^p F^{QH}_{t\frac{\partial}{\partial t}} \nonumber\\
&=t^p((\nabla^{QH}_{t\frac{\partial}{\partial t}})^p-\nabla^{QH}_{t\frac{\partial}{\partial t}}) \nonumber\\
&=t^p((\frac{1}{2}\mathrm{Deg}-\nabla^{QH}_{q\frac{\partial}{\partial q}})^p-(\frac{1}{2}\mathrm{Deg}-\nabla^{QH}_{q\frac{\partial}{\partial q}}))\nonumber\\
&=t^p((\frac{1}{2}\mathrm{Deg})^p-(\nabla^{QH}_{q\frac{\partial}{\partial q}})^p-(\frac{1}{2}\mathrm{Deg}-\nabla^{QH}_{q\frac{\partial}{\partial q}}))\qquad (\textrm{since $\nabla^{QH}_{q\frac{\partial}{\partial q}}$ has degree $0$, it commutes with $\mathrm{Deg}$})\nonumber\\
&=-F^{QH}_{tq\frac{\partial}{\partial q}},
\end{align}
where in the last inequality, $(\frac{1}{2}\mathrm{Deg})^p-\frac{1}{2}\mathrm{Deg}=0$ because $\mathrm{Deg}$ is multiplication by an integer.\qed

\begin{lemma}
For any $b\in H^*(X;\mathbf{k})$, $Q\Sigma_b$ commutes with $\nabla^{QH}_{t^2\frac{\partial}{\partial t}}$ of (1.2).     
\end{lemma}
\emph{Proof}. By Theorem 3.2, $Q\Sigma^q_b$ commutes with the quantum $q$-connection (3.22). Thus, by formula (B.2), it suffices to show that $Q\Sigma^q_b$ commutes with $\mathrm{Deg}$. However, $Q\Sigma^q_b$ is an operation of total degree $p|b|$. As a result, 
\begin{equation}
[\mathrm{Deg},Q\Sigma^q_b]=p|b|Q\Sigma^q_b=0    
\end{equation}
since we are in characteristic $p$. Restricting to $q=1$ gives the desired result.\qed

\renewcommand{\theequation}{C.\arabic{equation}}
\setcounter{equation}{0}
\section{An exponential type result for matrix factorizations}
Motivated by homological mirror symmetry \cite{Sh}, we consider the following $B$-side situation. Let $Y$ be an algebraic variety over $\overline{\mathbb{Q}}$ equipped with a function $W: Y\rightarrow \mathbb{A}^1$ with isolated singularities. % Fix $R\subset K$ a ring of integers of some number field (up to inverting finitely many elements) such that the pair $(Y,W)$ is defined.  
One can associate to the pair $(Y,W)$ its triangulated category of singularities \begin{equation}D^bSing(Y,W)=\prod_{\lambda\in\mathrm{crit(W)}}D^bSing(W^{-1}(\lambda)).
\end{equation}
The decomposition (C.1) should be thought of as the $B$-side analogue of the fact that there is one monotone Fukaya category $\mathrm{Fuk}(X)_{\lambda}$ associated with each eigenvalue $\lambda$ of $c_1\star$, cf. \cite[Section 2]{Sh}. By classical results of Orlov \cite{Or}, there is an exact equivalence of triangulated categories
\begin{equation}
D^bSing(W^{-1}(\lambda))\simeq D^b\mathrm{MF}(W-\lambda),
\end{equation}
where $\mathrm{MF}$ denotes the category of matrix factorizations. To simplify the computations, we restrict ourselves to considering a single summand of (C.1) and the following localized situation. Namely, let $Y=\mathrm{Spec}\;\overline{\mathbb{Q}}[z_1,\cdots,z_n]$ and $W:Y\rightarrow \mathbb{A}^1$ be a global function such that $W(0,\cdots,0)=dW(0,\cdots,0)=0$; moreover, assume that $(0,\cdots,0)$ is the only critical point of $W$.\par\indent
\cite[Theorem 1.1]{Shk} showed that the Getzler-Gauss-Manin $t$-connection on the periodic cyclic homology of $\mathrm{MF}(Y,W)$ is equivalent to the connection
\begin{equation}
\nabla^W_{\frac{\partial}{\partial t}}:=\frac{\partial}{\partial t}+\frac{W}{t^2}+\frac{\Gamma'}{t}    
\end{equation}
on $H^*(\Omega(Y)((t)),-dW+td)$, where $\Gamma'|_{\Omega^q(Y)}=-\frac{q}{2}$. Theorem C.1 below is an analogue of Theorem 1.3 in the context of matrix factorizations. This result is well known, see e.g. \cite[Theorem 1.1]{Sab}, but we give a different approach using reduction $p$ methods. 
\begin{thm}
The connection $\nabla^W_{\frac{\partial}{\partial t}}$ of (C.3) has regular singularity and quasi-unipotent monodromy at $t=0$. 
\end{thm}
\emph{Proof}. By Theorem A.1, it suffices to show that there exists a $R\subset \overline{\mathbb{Q}}$ finitely generated over $\mathbb{Z}$ over which $W$ is defined and such that when reduced mod each $\mathfrak{m}\in \mathrm{mSpec}(R)$, $\nabla^W_{\frac{\partial}{\partial t}}$ has nilpotent $p$-curvature (where $p$ equals the characteristic of the residue field $\kappa(\mathfrak{m})$). \par\indent
By Nullstellensatz, our assumption on $W$ implies that over $\overline{\mathbb{Q}}$ there exists a positive integer $N$ such that  
\begin{equation}W^N\in (\frac{\partial W}{\partial z_1},\cdots,\frac{\partial W}{\partial z_n}).\end{equation} 
Let $R$ be the ring obtained by adjoining to $\mathbb{Z}$ the coefficients of $W$ and of the polynomials appearing in expanding $W^N$ in terms of the generators $\{\frac{\partial W}{\partial z_1},\cdots,\frac{\partial W}{\partial z_n}\}$ as in $(C.4)$. We also adjoint $\frac{1}{2}$ to $R$, if it was not already contained. The argument goes in two steps.
\begin{enumerate}[label=\arabic*)]
    \item For $\mathfrak{m}\in\mathrm{mSpec}(R)$, we compute the $p$-curvature of $\nabla^W_{\frac{\partial}{\partial t}}$ over $\kappa(\mathfrak{m})$ along the vector field $t^2\frac{\partial}{\partial t}$. Because $W$ and $\Gamma'$ commute, 
     \begin{equation}
     F^W_{t^2\frac{\partial}{\partial t}}=(t^2\frac{\partial}{\partial t}+W+t\Gamma')^p
     =W^p+(t^2\frac{\partial}{\partial t}+t\Gamma')^p.
     \end{equation}
On the other hand,
\begin{align}
(t^2\frac{\partial}{\partial t}+t\Gamma')^p&=t^p((t\frac{\partial}{\partial t}+\Gamma')^p-(t\frac{\partial}{\partial t}+\Gamma'))\nonumber\\
&=t^p(t\frac{\partial}{\partial t}+(\Gamma')^p-t\frac{\partial}{\partial t}-\Gamma')\nonumber\\
&=0,
\end{align}
where in the last equality, we used that $(t\frac{\partial}{\partial t})^p=t\frac{\partial}{\partial t}$ and $(\Gamma')^p=\Gamma'$ (since $\Gamma'$ always takes half-integer values) in characteristic $p>2$. 
\item By 1), it suffices to show that multiplication by $W^p$ defines a nilpotent operation on $H^*(\Omega(Y)((t)),-dW +td)$
. The key observation is that there is a Frobenius $p$-linear algebra action of the Jacobian ring $[f]\in\kappa(\mathfrak{m})[z_1,\cdots,z_n]/(\partial_{z_1} W,\cdots,\partial_{z_n}W)$ on `twisted de Rham cohomology' $[\alpha]\in H^*(\Omega(\kappa(\mathfrak{m})[z_1,\cdots,z_n])((t)),-dW+td)$ given by
\begin{equation}
\mathrm{Act}_{[f]}^p([\alpha]):=[f^p\cdot \alpha].
\end{equation}
We check this is indeed well-defined. $d(f^p)=0$ implies that $f^p$ commutes with $-dW+td$. On the other hand, if $f=\iota_{dW}(D)=D(W)$ for some vector field $D$, a straightforward but tedious computation shows that
$-D(W)^p\alpha+t^{p-1}D^p(W)\alpha+t^{p}\mathcal{L}_{D^p}\alpha=$
\begin{equation}
[-dW+td, \sum_{\substack{\sum_{j=1}^se_j=
\sum_{j=1}^se_jm_j+l=p-1\\
e_j>0,\,l\geq 0, \,0\leq m_1<m_2<\cdots<m_s}}\frac{(p-1)!^2 e(m_1)}{l!\prod_{j=1}^s(m_j!)^{e_j}(e_j!)} D^{m_1}(W)^{e_1}\cdots D^{m_s}(W)^{e_s}\mathcal{L}_D^{l}\iota_D] \alpha,
\end{equation}
where we artificially define $D^0(W):=-t$, and $e(m_1):=e_1!$ if $m_1=0$ and $1$ if $m\neq 0$. If $D\in\{\frac{\partial}{\partial z_1},\cdots,\frac{\partial}{\partial z_n}\}$, then $D^p=0$, and thus (C.8) implies that $D^p(W)\alpha$ is a coboundary if $(-dW+td)\alpha=0$; since the coordinate vector fields generate, this holds for all $D$. \par\indent
As a consequence, in order to show that $\mathrm{Act}^p_W=W^p\cdot$ defines a nilpotent operation, it suffices to show that $W\in \kappa(\mathfrak{m})[z_1,\cdots,z_n]/(\partial_{z_1} W,\cdots,\partial_{z_n}W)$ is a nilpotent element, which follows from (C.4).
\end{enumerate}\qed

\end{appendices}

\end{document}